\documentclass[review,12pt]{elsarticle}
\usepackage{amsmath}
\usepackage{amsfonts}
\usepackage{geometry}   
\usepackage[parfill]{parskip}    
\usepackage{graphicx}
\usepackage{amssymb}
\usepackage{epstopdf}
\usepackage{psfrag}
\usepackage{natbib}

\newcommand {\bea}{\begin{eqnarray}}
\newcommand {\ea}{\end{eqnarray}}
\newtheorem{thm}{Theorem}[section]
\newtheorem{prob}{Problem}[section]

\newenvironment{proof}[1][Proof]{\textbf{#1.} }{\hspace{\stretch{1}}\rule{0.5em}{0.5em}}

\usepackage{xspace}

\usepackage{subfigure}
\usepackage[all]{xy}
\newcommand{\thmref}[1]{{Theorem~\ref{#1}}}

\usepackage{subfig}
 \usepackage{float}
\usepackage{latexsym}
\usepackage{xcolor}
\usepackage{amsfonts}
\usepackage{graphicx}
\usepackage{cases}
\usepackage{multirow}
\usepackage{amssymb}

\usepackage{fancyhdr}
\journal{Mathematics and Computers in Simulation}
\begin{document}
\begin{frontmatter}
\title{A fitted finite volume method for stochastic optimal control Problems}

\author[jdm]{Christelle Dleuna Nyoumbi}
\ead{christelle.dleuna@imsp-uac.org, christelle.d.nyoumbi@aims-senegal.org}
\address[jdm]{Institut de Math\'{e}matiques et de Sciences Physiques (IMSP), Universit\'{e} d'Abomey-Calavi\\ 01 B.P. 613, Porto-Novo, Benin.}

\author[at,atb,atc]{Antoine Tambue}
\cortext[cor1]{Corresponding author}
\ead{antonio@aims.ac.za}
\address[at]{Department of Computer science, Electrical engineering and Mathematical sciences, Western Norway University of Applied Sciences, Inndalsveien 28, 5063 Bergen.}
\address[atb]{Center for Research in Computational and Applied Mechanics (CERECAM), and Department of Mathematics and Applied Mathematics, University of Cape Town, 7701 Rondebosch, South Africa.}
\address[atc]{The African Institute for Mathematical Sciences(AIMS),
6-8 Melrose Road, Muizenberg 7945, South Africa.}

%

\begin{abstract}
In this article, we provide a numerical method  based on fitted finite volume method to approximate the Hamilton-Jacobi-Bellman (HJB) equation
    coming from stochastic optimal control problems.  The computational challenge  is due to the nature of the
    HJB equation, which  may be a second-order degenerated partial differential equation  coupled with optimization. 
    In the work,  we discretize the HJB equation  using the fitted finite volume  method and show that  matrix  resulting  from spatial discretization  is an M-matrix. 
    The optimization  problem  is solved  at every time step using iterative method.
     Numerical results are presented to show  the robustness of the fitted finite volume  numerical method comparing to the standard finite difference method.\\
\end{abstract}
\end{frontmatter}
\section{Introduction}
We consider  the numerical approximation of the following  controlled Stochastic Differential Equation (SDE) defined  in $\mathbb{R}^{ n}$  ($ n\geq 1$)  by
     \begin{align} \label{premier}
    dx_t = b(t,x_t, \alpha_t) dt + \sigma (t, x_t, \alpha_t) d\omega_t,\,\, x(0)=x_0 
\end{align}
where  
\begin{align}
	b&: [0, T] \times \mathbb{R}^n \times \mathcal{A} \rightarrow  \mathbb{R}^n \nonumber\\ 
	& (t,x_t, \alpha_t)) \rightarrow b(t,x_t,\alpha_t)
\end{align}
is the drift term  and
\begin{align}
	\sigma &:[0, T] \times \mathbb{R}^n\times \mathcal{A} \rightarrow  \mathbb{R}^{n\times d} \nonumber \\
	& (t,x_t,\alpha_t) ) \rightarrow \sigma(t,x_t,\alpha_t)
\end{align}
the  d-dimensional diffusion coefficients.  Note that $ \omega_t $ are d-dimensional independent Brownian motion on $ \left(\Omega, \mathcal{F}, (\mathcal{F}_t)_{t \geq 0}, \mathbb{P} \right) $, $ \alpha =(\alpha_t)_{t \geq 0}$ is an $ \mathbb{F} $-adapted process, valued in  $ \mathcal{A} $ closed convex subset of $ \mathbb{R}^m \, (m \geq 1)$ and satisfying some integrability conditions and/or state constraints. Precise assumptions on $ b $ and $ \sigma $  to ensure the existence of the unique solution $
  x_t $ of (\ref{premier})  can be found in \cite{HP}.\\
Given a function $ g $ from  $\mathbb{R}^n $ into  $ \mathbb{R}$ and $ f $ from  $[0, T]\times \mathbb{R}^n \times \mathcal{A} $ into  $ \mathbb{R}$,  the value function is defined by 
 \begin{align} \label{vienu}
	v(t, x) = \underset{\alpha \in \mathcal{A}}{\sup}\, \mathbb{E}\, \left[ \int_{t}^T f(s, x, \alpha) \,ds +g (x_T)\right], \,\,\,x \in \mathbb{R}^n, 
\end{align} 
and the resulting Hamilton Jacobi-Bellamn (HJB) equation (see \cite{KL})is given by
\begin{align}\label{merci}
	\begin{cases}
		v_t(t, x) + \underset{\alpha \in \mathcal{A}}{\sup} \left[L^{\alpha} v(t, x) + f(t, x,\alpha)\right] = 0 \quad\text{on} \  [0,T)\times \mathbb{R}^n\\
		v(T,x) = g(x), \,\,\,\,x \,\in \mathbb{R}^n
	\end{cases}
\end{align}
where 
\begin{align}\label{tous}
	L^{\alpha} v(t, x) = \sum_{i=1}^n (b(t, x,\alpha))_i \dfrac{\partial v(t,x)}{\partial x_i}  + \sum_{i,j=1}^n ( a (t,x,\alpha))_{i,j}\,\dfrac{\partial^2 v(t,x)}{\partial x_i\,\partial x_j},
\end{align}   
and $ a(t,x,\alpha) = (\dfrac{1}{2}(\sigma (t,x,\alpha))(\sigma(t,x,\alpha))^T )_{i,j} $.\\
The existence and uniqueness of the  viscosity  solution  of the HJB equation (\ref{merci})  is well  known and  can be found in \cite{HP}.
Equation (\ref{merci}) is a  initial value problem. There are two
unknown functions in this equation,  the value function $ v $ and the optimal control $ \alpha $. 
In most practical situations, (\ref{merci}) is not analytically solvable therefore numerical approximations are 
the only tools appropriate to provide reasonable approximations. 
Numerical approximation of HJB-equation of type (\ref{merci}) is therefore an active research area and has attracted a lot of attentions \cite{Cra7,KNV1,KL1,KL, Cra2, Cra1, Cra3,JER,HPFH}. 
While  solving numerically  HJB equation, the keys challenge are the low regularity of the solution of HJB equation and the lack of appropriate  numerical 
methods to tackle the degeneracy of the differential operator in HJB equation. Indeed adding to the standard issue that we usually have when solving degenerated PDE, we need to couple 
with  an optimization problem at every grid  point and every time step. In terms of existing numerical methods, there are two basic threads of literature concerning
controlled HJB equations. A standard approach is based on  Markov chain approximation. In financial terms, this approach is equivalent to an explicit finite difference method. 
However, these methods are well-known to suffer from time step limitations due to stability issues \cite{Peter}. A more recent approach is based on numerical methods  such as 
 finite difference method which ensure convergence to the viscosity solution of the HJB equation couple with an optimization problem at each time \cite{IGA} .  

For many stochastic optimal control  problems such as Merton's control problem, the linear operator is degenerated  when the spatial variables approach
the region near to zero.
This degeneracy has an adverse impact on the accuracy when the finite difference method is used to solve the PDE (see \cite{wilmott2005best}, chapter 26).
This  degeneracy  also has an adverse impact on the accuracy of our stochastic optimal control  problems since  its numerical resolution  implies  the resolution of PDE,  coupled with optimization problem.

In this article, we propose a numerical scheme based on a finite volume method  suitable to handle the degeneracy of the linear operator  while solving  numerically  the HJB equation in dimension $1$ and $2$.
The method  is coupled with  implicit time-stepping method for temporal discretization  method and the iterative method  presented in \cite{HPFH} for optimization  problem at every time step. 
 More precisely, this method is based on  fitted  finite volume technique   proposed in \cite{huangfitted2006} to solve the degenerated  Black Sholes equations.  
 Note that to the best of our knowledge, such method  has not been used to solve the stochastic optimal control problem \eqref{merci}.
 
 The merit of the method  is  that it is absolutely stable in time because of the implicit nature of the time discretisation and  the corresponding matrix after spatial discretization  is  a positive-definite $M $-matrix.\\
 Numerical simulations prove  that our proposed method is more accurate that  the standard method based on finite difference spatial discretization.

The rest of this article is organized as follows.  In section \ref{sec1}, we  present the finite volume method with the fitting technique for dimension 1 and 2.  We will also show that the system matrix of the resulting
discrete equations is an $M $-matrix. In section \ref{sec2}, we will present the temporal discretization and optimization problem in dimension 1 and 2. Numerical experiments using Matlab software will be performed in section \ref{sec3} 
to demonstrate the accuracy of the proposed numerical  method. We conclude the work at section \ref{sec4}  by summarizing our finding.
  
  \section{Spatial discretization}
  \label{sec1}
  As we already know,  the resolution  of the  HJB  equation \eqref{merci}  involves   a spatial discretisation, a temporal discretisation and  an optimisation problem  at every  grid  point and  each time step.
  The goal of this section is to provide the spatial discretization of the HJB equation \eqref{merci} solving our stochastic optimal control problem \eqref{vienu}. 
   Details  in this section can be found  in \cite{WS}, where such methods have been used to solve  the degenerated Black Sholes equation for option pricing with constant coefficients.
  \subsection{  Spatial discretization   based on  fitted  finite volume method in dimension 1}
Consider the more generalized HJB equation \eqref{merci} in  dimension 1 ($n=1$) which  can be written in the form.  
\begin{align} \label{vie1}
	\dfrac{\partial v(x,t)}{\partial t} + \sup_{\alpha \in \mathcal{A}}\left[ \dfrac{\partial }{\partial x} \left(a (x,t,\alpha)\,x^2 \dfrac{\partial v(x,t)}{\partial x}+ b (x,t,\alpha)\,x\,v(x,t) \right)+c(x,t,\alpha)\,v(x,t)\right] = 0,
\end{align}
where $a(t,x,\alpha)>0$, $\alpha=\alpha(x,t)$ and bounded.
As usual, we truncate the problem  in the finite interval  $I= [0,x_{\text{max}}] $.  Let the interval $ I= [0,x_{\text{max}}] $ be divided into $ N_1 $ sub-intervals
$  I_i :=(x_i,x_{i+1}),\,\,\,\,\, i=0\cdots N_1-1$
with $ 0 =x_0 < x_1< \cdots \cdots< x_{N_1}= x_\text{max} $. We also set $ x_{i+1/2} =\dfrac{x_{i} +x_{i+1} }{2} $ and $ x_{i-1/2} =\dfrac{x_{i-1} +x_{i} }{2} $ 
for each $ i=1\cdots N_1-1$ . If we define $ x_{-1/2} = x_0$ and $ x_{N_1+1/2} = x_\text{max}$ integrating both size of (\ref{vie1}) over $ J_i=\left( x_{i-1/2},\, x_{i+1/2}\right) $  and  taking $\alpha_i=\alpha(x_i,t)$,  we have 
\begin{align}
	 \int_{x_{i-1/2}}^{x_{i+1/2}}\dfrac{\partial v }{\partial t} dx + \int_{x_{i-1/2}}^{x_{i+1/2}} \sup_{\alpha_i \in \mathcal{A}}\left[ \dfrac{\partial }{\partial x}  x \left(a(x,t,\alpha_i)\,x \dfrac{\partial v}{\partial x}+ b(x,t,\alpha_i) \,v \right)+c(x,t,\alpha_i)\,v\right]dx = 0
\end{align}

Applying the mid-points quadrature rule to the first and the last point terms, we obtain the above
\begin{align}\label{bjrx}
	\dfrac{d v}{d t} l_i + \sup_{\alpha_i \in \mathcal{A}}\left[\left[x_{i+1/2}\rho(v)\left|_{x_{i+1/2}}-x_{i-1/2} \rho(v)\right|_{x_{i-1/2}}\right]+c(x_i,t,\alpha_i)\, v_i\, l_i \right] =0,
\end{align}
for $ i= 1,2,\cdots N_1-1 $, where $ l_i = x_{i+1/2} - x_{i-1/2} $ is the length of $ J_i $. $ v_i $ denotes the nodal approximation to $ v(\tau, x_i) $ and $ \rho(v) $ is the flux associated with $ v $ defined by
\begin{align}
	\rho(v) := a(x,t,\alpha_i)\,x\,\dfrac{\partial v}{\partial x} + b(x,t,\alpha_i)\,v.
\end{align}
Clearly, we now need to derive approximation of the flux defined above at the mid-point $ x_{i+1/2} $, of the interval $ I_i $ for $ i= 2,\cdots N_1-1.$ 
This discussion is divided into two cases for $ i \geq 1 $ and $ I_0 = (0,x_1) $.\\
\textbf{\underline{Case I}:} Approximation of $ \rho $ at $ x_{i+1/2} $ for $ i\geq 2 $.\\
The term $ \left(a (x,t,\alpha_i)\,x \dfrac{\partial v}{\partial x}+ b(x,t,\alpha_i)\,v \right) $  is approximated by solving the  boundary value problem
\begin{align}\label{vie2x}
	\left(a(x,t,\alpha_i)\,x \dfrac{\partial v}{\partial x}+ b(x_{i+1/2},t,\alpha_i),v, \right)'= 0,\, \,\,\, x \in I_i\\
	v(x_i)= v_i(t),\,\, v(x_{i+1})= v_{i+1}(t).
\end{align}
Integrating (\ref{vie2x}) yields the first-order linear equations
\begin{align}
	\rho_i(v)(t)= a(x,t,\alpha_i)\,x \dfrac{\partial v}{\partial x}+ b(x_{i+1/2},t,\alpha_i)\,v = C_1
\end{align}
where $ C_1 $ denotes an additive constant.  As in \cite{WS}, the solution is given  by

\begin{align}
	v (t)= \dfrac{C_1}{b(x_{i+1/2},t,\alpha_i)} + C_2\, x^{-\dfrac{b(x_{i+1/2},t,\alpha_i)}{a(x_{i+1/2},t,\alpha_i)}}.
\end{align}
Note that in this deduction we have assumed that $ b(x_{i+1/2},t,\alpha_i) \neq 0 $.   By setting  $ \beta_i(t) =\dfrac{b(x_{i+1/2},t,\alpha_i)}{a(x_{i+1/2},t,\alpha_i)} $,  using the boundary conditions in (\ref{vie2x}) yields
\begin{align}
	v_i (t)= \dfrac{C_1}{b(x_{i+1/2},t,\alpha_i)} + C_2\, x_i^{-\beta_i(t)}\,\,\,	\text{and} \,\,\,v_{i+1} (t)= \dfrac{C_1}{b(x_{i+1/2},t,\alpha_i)} + C_2\, x_{i+1}^{-\beta_i(t)}
\end{align}
Solving the   following linear system  with respect to $C_1$ and $C_2$ yields
\begin{align}
	\begin{cases}
		v_i (t)= \dfrac{C_1}{b(x_{i+1/2},t,\alpha_i)} + C_2\, x_i^{-\beta_i(t)} \\
		v_{i+1}(t) = \dfrac{C_1}{b(x_{i+1/2},t,\alpha_i)} + C_2\, x_{i+1}^{-\beta_i(t)}
	\end{cases}
\end{align}
yields
\begin{align}
	\rho_i(v)(t) = C_1 = \dfrac{b(x_{i+1/2},t,\alpha_i)\,\left(x_{i+1}^{\beta_i(t)}v_{i+1}(t)-x_{i}^{\beta_i(t)}v_i (t)\right)}{x_{i+1}^{\beta_i(t)}-x_{i}^{\beta_i(t)}}
\end{align}
$\rho_i (v)(t) $ provides an approximation to the $ \rho(v)(t) $ at $ x_{i+1/2} $.\\
\textbf{\underline{Case II:}}   This is the degenerated zone.   The aims  here  is  to  approximate  $ \rho $ at $ x_{1/2} $  in the  sub-interval $ I_0 $.
In this case,   the  following  problem  is considered
\begin{align}\label{vie4x}
	(a(x_{1/2},t,\alpha_1)\,x \dfrac{\partial v}{\partial x}+b(x_{1/2},t,\alpha_1)\,v )'= C_2\,\,\,\textbf{in}\,\,[0,x_1] \\ \nonumber
	v(0) = v_0(t),\,\,\,v(x_1)= v_1(t)
\end{align}
where $ C_2 $ is an unknown constant to be determined.  Following \cite{WS}, integrating (\ref{vie4x}) yields
\begin{eqnarray}
	\rho_0(v)\vert_{1/2} (t)= a(x_{1/2},t,\alpha_i)\,x_{1/2} \dfrac{\partial v}{\partial x}+ b(x_{1/2},t,\alpha_1)\,v = b(x_{1/2},t,\alpha_1)\,v_0(t) +C_2\,x_{1/2}.
\end{eqnarray}
Since $ x_{1/2} = \dfrac{x_1+ x_0}{2}$  with  $x_0=0$, we have   $ C_2\,x_1 = (a(x_{1/2},t,\alpha_1)+b(x_{1/2},t,\alpha_1))(v_1(t)-v_0(t))$.
Therefore we have  
\begin{align}
	\rho_0(v)\vert_{1/2} (t)=  \dfrac{1}{2}\left[ (a(x_{1/2},t,\alpha_1)+ b(x_{1/2},t,\alpha_1))v_1(t)-(a(x_{1/2},t,\alpha_1)-b(x_{1/2},t,\alpha_1))v_0(t)\right].
\end{align}
By replacing $ \rho  $ by its approximated value, (\ref{bjrx}) becomes for $ i = 0,1,\cdots,N_1-1 $
{\small{
\begin{align}
	\dfrac{d v_i(t) }{ dt}  + \sup_{\alpha_i \in \mathcal{A}} \dfrac{1}{l_i}\,\left[x_{i+1/2}\dfrac{b(x_{i+1/2},t,\alpha_i) \,\left(x_{i+1}^{\beta_i(t)}v_{i+1}(t)-x_{i}^{\beta_i(t)}v_i(t) \right)}{x_{i+1}^{\beta_i(t)}-x_{i}^{\beta_i(t)}}\right.\\ \nonumber
	\left. -x_{i-1/2}\dfrac{b(x_{i-1/2},t,\alpha_i)\,\left(x_{i}^{\beta_{i-1}(t)}v_{i}(t)-x_{i-1}^{\beta_{i-1}(t)}v_{i-1} (t)\right)}{x_{i}^{\beta_{i-1}(t)}-x_{i-1}^{\beta_{i-1}(t)}}+c_i(t,\alpha_i)\, v_i (t)\,l_i \right] = 0
\end{align}
}}
By setting $ \tau=T-t $ and including the boundary conditions, we have  the following  
system of Ordinary Differential Equation (ODE) coupled with optimisation  problem. 
\begin{align}
\label{p4xx}
	\begin{cases}
		- v_{\tau}(\tau) + \underset{\alpha \in \mathcal{A}^{N_{1}-1}}{\sup}\left[A(\alpha,\tau)\,v(\tau) + G( \alpha, \tau) \right]= 0\\
		v(0) \;\;\text{given},
	\end{cases}
\end{align}
which can be rewritten as   
\begin{align}\label{p4x}
	\begin{cases}
		 v_{\tau}(\tau) + \underset{\alpha \in  \mathcal{A}^{N_{1}-1}}{\inf}\left[E(\alpha,\tau)\,v(\tau)+ F( \alpha, \tau) \right] = 0\\
		v(0) \;\;\text{given},
	\end{cases}
\end{align}
where  $v(\tau)=(v_1(\tau),\cdots, v_{N_1-1}(\tau)) $ and $F(\alpha,\tau)= (F_1(\alpha_1,\tau),\cdots, F_{N_1-1}(\alpha_{N_1-1},\tau)) $ includes all  Dirichlet boundary and final conditions,  $A (\alpha,\tau) = -E(\alpha,\tau) $  and $ G(\alpha,\tau) = - F(\alpha,\tau)$  are defined as for  $i = 1,\cdots,N_1-1 $\\
\vspace*{0.3cm}
\begin{align}
	E_{i,i+1}(\alpha_i, \tau)& = -x_{i+1/2}\dfrac{b_{i+1/2}(\tau,\alpha_i)\,x_{i+1}^{\beta_i(\tau)} }{l_i\,(x_{i+1}^{\beta_i(\tau)}-x_{i}^{\beta_i(\tau)})},\\
	E_{i,i}(\alpha_i, \tau) &= \left(x_{i+1/2}\dfrac{b_{i+1/2}(\tau,\alpha_i)\,x_{i}^{\beta_i(\tau)} }{l_i\,(x_{i+1}^{\beta_i(\tau)}-x_{i}^{\beta_i(\tau)})}+ x_{i-1/2}\dfrac{b_{i-1/2}(\tau,\alpha_i)\,x_{i}^{\beta_{i-1}(\tau)} }{l_i\,(x_{i}^{\beta_{i-1}(\tau)}-x_{i-1}^{\beta_{i-1}(\tau)})}- c_i(\tau,\alpha_i) \right),\\
	E_{i,i-1}(\alpha_i, \tau)& = -x_{i-1/2}\dfrac{b_{i-1/2}(\tau,\alpha_i)\,x_{i-1}^{\beta_{i-1}(\tau)} }{l_i\,(x_{i}^{\beta_{i-1}(\tau)}-x_{i-1}^{\beta_{i-1}(\tau)})},
\end{align}
\begin{align}
	E_{1,1}(\alpha_1, \tau)& =   x_{1+1/2}\dfrac{b_{1+1/2}(\tau,\alpha_1)\,x_{1}^{\beta_1(\tau)} }{l_1\,(x_{2}^{\beta_1(\tau)}-x_{1}^{\beta_1(\tau)})} + \dfrac{1}{4\,l_1} x_{1}(a_{1/2}(\tau,\alpha_1)+b_{1/2(\tau,\alpha_1)}) - c_1(\tau,\alpha_1) \\
	E_{1,2}(\alpha_1, \tau)& = - x_{1+1/2}\dfrac{b_{1+1/2}(\tau,\alpha_1)\,x_{2}^{\beta} }{l_1\,(x_{2}^{\beta_1(\tau)}-x_{1}^{\beta_1(\tau)})}
\end{align}

\vspace*{0.3cm}
$ G( \alpha, \tau)= \left[\begin{array}{c}
-\dfrac{1}{4\,l_1}\,x_{1}(a_{1/2}(\tau,\alpha_1)-b_{1/2}(\tau,\alpha_1))\,v_0 \\ 
0\\ 
\vdots\\ 
0\\ 
-x_{N_1-1/2}\dfrac{b_{N_1-1/2}(\tau,\alpha_{N_1-1})\,x_{N_1}^{\beta_i(\tau)} }{l_{N_1-1}\,(x_{N_1}^{\beta_{N_1-1}(\tau)}-x_{N_1-1}^{\beta_{N_1-1}(\tau)})}v_{N_1}
\end{array} \right].$ \\

	\begin{thm} \label{tm} 
		      Assume that $ c_i(\tau,\alpha) < 0 $, $i = 1, \cdots ,N_1-1 $,  let  $h=\underset{1\leq i\leq N_1} {\max} l_i$. If $h$  is relatively small then the matrix $E (\alpha, \tau) $ 
		      in the system (\ref{p4}) is an \textbf{M}-matrix for any $ \alpha \,\in\,\mathcal{A}$. 
	\end{thm}
	\hspace*{0.5cm}\\
	\begin{proof} 
	Let us show that $E(\alpha, \tau) $ has positive diagonal, non-positive off diagonal, and is diagonally dominant. We first note that 
		\begin{align}
			\dfrac{b_{i+1/2}(\tau,\alpha)}{x_{i+1}^{\beta_i(\tau)}- x_{i}^{\beta_i(\tau)}} = \dfrac{a_{i+1/2}(\tau,\alpha)\,\beta_i(\tau) }{x_{i+1}^{\beta_i(\tau)}- x_{i}^{\beta_i(\tau)}} > 0,
		\end{align}
		for $ i = 1,\cdots, N_1-1 $,  and all $ b_{i+1/2}(\tau,\alpha) \neq 0, $ $ b_{i-1/2}(\tau,\alpha) \neq 0, $\,\, with $ a_{i+1/2}(\tau,\alpha)> 0$ and $ a_{i-1/2}(\tau,\alpha)> 0 $.
	
	This also holds when $ b_{i+1/2}(\tau,\alpha) \rightarrow 0 $ and $ b_{i-1/2}(\tau,\alpha) \rightarrow 0 $, that is 
		\begin{align}
	&\lim_{b_{i+1/2}(\tau,\alpha)\rightarrow 0}	\dfrac{b_{i+1/2}(\tau,\alpha)}{x_{i+1}^{\beta_i(\tau)}- x_{i}^{\beta_i(\tau)}} = \dfrac{b_{i+1/2}(\tau,\alpha) }{e^{\beta_{i}(\tau)\ln(x_{i+1})}-e^{\beta_{i}(\tau)\ln(x_{i})}}= \dfrac{b_{i+1/2}(\tau,\alpha) }{\beta_{i}(\tau)\ln(x_{i+1})-\beta_{i}(\tau)\ln(x_{i})} \\ \nonumber
	&= a_{i+1/2}(\tau,\alpha)\left(\ln \dfrac{ x_{i+1}}{x_{i}}\right)^{-1} > 0,\\ \nonumber  & \,\lim_{b_{i-1/2}(\tau,\alpha)\rightarrow 0}	\dfrac{b_{i-1/2}(\tau)}{x_{i}^{\beta_{i-1}(\tau)}- x_{i-1}^{\beta_{i-1}(\tau)}} = \dfrac{b_{i-1/2}(\tau,\alpha) }{e^{\beta_{i-1}(\tau)\ln(x_{i})}-e^{\beta_{i-1}(\tau)\ln(x_{i-1})}}= \dfrac{b_{i-1/2}(\tau,\alpha) }{\beta_{i-1}(\tau)\ln(x_{i})-\beta_{i-1}(\tau)\ln(x_{i-1})}\\ \nonumber
	&= a_{i-1/2}(\tau,\alpha) \left(\ln \dfrac{ x_{i}}{x_{i-1}}\right)^{-1} > 0 
		\end{align}.
		Using the definition of  $ E (\alpha,\tau) $ given above, we see that 
		\begin{align*}
			E_{i,i} \geqslant 0,\,\,\,E_{i,i+1} \leqslant 0,\,\,\, E_{i,i-1}\leqslant 0\;\;\;  i = 2,\cdots,N_1-1,
		\end{align*}
		\begin{align*}
			\left|E_{i,i}\right| & \geq \left|E_{i,i-1}\right|+\left|E_{i,i+1}\right|
		\end{align*}
		because $x_{i+1}^{\beta_i(\tau)} \approx x_{i}^{\beta_i(\tau)} + x_{i}^{\beta_i(\tau)-1}\,{\beta_i(\tau)}\,h$,\,\, $x_{i-1}^{\beta_{i-1}(\tau)} \approx x_{i}^{\beta_{i-1}(\tau)} - x_{i}^{\beta_{i-1}(\tau)-1}\,{\beta_{i-1}(\tau)}\,h$ \,\,and
		\begin{align*}
		& \left|E_{i,i}\right| - \left|E_{i,i+1}\right|-\left|E_{i,i-1}\right|\\ & = -\underset{\rightarrow 0}{\underbrace{\underset{>0}{\underbrace{\left(\dfrac{b_{i+1/2}(\tau)}{x_{i+1}^{\beta_i(\tau)}-x_{i}^{\beta_i(\tau)}}\right)}} \underset{\rightarrow 0}{\underbrace{\left(h^{\beta_i}\,\beta_i\,x_{i}^{\beta_{i}-1}\right)}} }}+ \underset{\rightarrow 0}{\underbrace{\underset{>0}{\underbrace{\left(\dfrac{b_{i-1/2}(\tau,\alpha)}{x_{i}^{\beta_{i-1}(\tau)}-x_{i-1}^{\beta_{i-1}(\tau)}}\right)}} \underset{\rightarrow 0}{\underbrace{\left(h^{\beta_{i-1}}\,\beta_{i-1}\,x_{i}^{\beta_{i-1}-1}\right)}} }}  - c_i(\tau,\alpha).
		\end{align*}
		 Note that for $ i = 1$,  we  have $ E_{1,1}\geq 0 $\,\,if\, ${a_{1/2}(\tau,\alpha)} + b_{1/2}(\tau,\alpha) $,\,\,are nonnegative and $ c_1(\tau,\alpha) < 0 $. 
		 So  $ E(\alpha, \tau)$ is diagonally dominant  and  is therefore  an  $\textbf{M} $-matrix.
\end{proof}
 \subsection{ Spatial discretization based on fitted finite volume method in dimension 2}

    Here we consider the following two dimensional problem
  
    	\begin{align} \label{vie1y}
    		\dfrac{\partial v(x,y,t) }{\partial t} + \sup_{\alpha \in \mathcal{A}}\left[ \nabla\cdot \left( k(x,y,t,\alpha) \right) + c(x,y,t,\alpha)\,v(x,y,t) \right] = 0,
    	\end{align}
    	where $k (x,y,t,\alpha)= A(x,y,t,\alpha)\cdot\nabla v(x,y,t)+ b\,v(x,y,t)$ is the flux, 
	
	$b = (x\,b_1(x,y,t,\alpha), y\, b_2(x,y,t,\alpha))^T$ and
    	\[
    	A=\left( \begin{array}{cc}
    	a_{11} & a_{12} \\
    	a_{21} & a_{22}
    	\end{array} \right),
    	\]
 We will  assume  that $a_{21}= a_{12}$.  We also  also define the following coefficients, which  will help us to build our scheme smoothly 
$ a_{11}(x,y,t,\alpha) = a(x, y, t, \alpha)\,x^2, a_{22}(x,y,t,\alpha) = \overline{a}(x,y,t,\alpha) y^2$ and $  a_{1 2}=a_{2 1} =d_1(x,y,t,\alpha) xy$.
	
As usual the two dimensional domain is truncated  to $I_{x} \times I_{y}$, where  $ I_{x} = [0,x_{\text{max}}] $ and $ I_{y}= [0,y_{\text{max}}]  $ be divided into $ N_1 $ and $ N_2 $ sub-intervals:
    	$$  I_{x_{i}} :=(x_{i}, x_{i+1}),\,\, I_{y_j} := (y_{j}, y_{j+1}),\,\,\,\, i=0,\cdots ,N_1-1,\,\,j=0,\cdots, N_2-1$$
    	with $ 0 = x_{0} < x_{1}< \cdots \cdots< x_{N_1}= x_{\text{max}} $ and $ 0 = y_{0} < y_1 < \cdots \cdots< y_{N_2}= y_{\text{max}} $. 
	This defines a mesh on $ I_{x} \times I_{y} $ with all the mesh lines perpendicular to one of the
    	axes. \\
  We also set  $$  x_{i+1/2} =\dfrac{x_{i} + x_{i+1} }{2},\, x_{i-1/2} =\dfrac{x_{i-1} + x_{i} }{2},\, y_{j+1/2} =\dfrac{y_{j} + y_{j+1} }{2},\, \; y_{j-1/2} =\dfrac{y_{j-1} + y_{j} }{2}, $$ for each $ i=1,\cdots, N_1-1$\,and each $ j=1,\cdots, N_2-1$. We denote $ N = (N_1-1)\times (N_2-1)$. These mid-points form a second partition of $ I_{x} \times I_{y} $ if we define $ x_{-1/2} = x_{0}$, $ x_{N_1+1/2} = x_{\text{max}}$,\, $ y_{-1/2} = y_{0}$ and $ y_{N_2+1/2} = y_{\text{max}}$. For each $i = 0, 1, \cdots ,N_1 $ and $j = 0, 1,\cdots,N_2 $, we set $h_{x_i} = x_{i+1/2} - x_{i-1/2} $ and  $h_{y_j}  = y_{j+1/2} - y_{j-1/2}$. \\
    We now discuss the finite volume method for (\ref{vie1y}). Integrating both size of (\ref{vie1y}) over $ \mathcal{R}_{i,j}=\left[ x_{i-1/2},x_{i+1/2} \right] \times \left[  y_{j-1/2},y_{j+1/2}\right] $, we have
    \begin{align}
    		\int_{x_{i-1/2}}^{x_{i+1/2}} \int_{y_{j-1/2}}^{y_{j+1/2}}\dfrac{\partial v }{\partial  t}\, dx\, dy + \int_{x_{i-1/2}}^{x_{i+1/2}}\int_{y_{j-1/2}}^{y_{j+1/2}} \sup_{\alpha \in \mathcal{A}} \left[ \nabla \cdot \left( k(v (x,y,t,\alpha))\right) + c(x,y,t,\alpha)\,v(x,y,t) \right]\, dx\,dy = 0,
    	\end{align}
    	for $ i =1,2,\cdots, N_1-1 $,\, $ j =1,2,\cdots, N_2-1 $.
    	Applying the mid-points quadrature rule to the first and the last point terms, we obtain the above
    	\begin{align}\label{bjr}
    		\dfrac{d v_{i,j} (t) }{d t}\,l_{i,j} + \sup_{\alpha_{i,j} \in \mathcal{A}} \left[ \int_{\partial \mathcal{R}_{i,j}}\nabla \cdot \left( k(v(x,y,t,\alpha_{i,j} )\right)\,dx\,dy + c_{i,j}(t,\alpha_{i,j} )\, v_{i,j}(t)\,l_{i,j}\right] = 0
    	\end{align}
    	for $ i= 1,2,\cdots N_1-1 $,\, $ j= 1,2,\cdots N_2-1 $ where $ l_{i,j} = \left( x_{i+1/2}-x_{i-1/2} \right) \times \left(  y_{j+1/2}-y_{j-1/2}\right) $ is the length of $ \mathcal{R}_{i,j} $, and $ v_{i,j}(t) $ denotes the nodal approximation to $ v( x_{i}, y_{j},t) $.
    		We now consider the approximation of the middle term in (\ref{bjr}). Let $ \bf n $ denote the unit vector outward-normal to $ \partial \mathcal{R}_{i,j} $. By Ostrogradski Theorem, integrating by parts and using the definition of flux $ k $, we have
    	\begin{align}\label{vol}
    		\int_{\mathcal{R}_{i,j}} \nabla \cdot \left( k(v)\right) & = \int_{\partial \mathcal{R}_{i,j}} k(v(x,y,t,\alpha_{i,j})) \cdot \bf n \, ds \nonumber \nonumber\\
    		&= \int_{\left(x_{i+1/2},y_{j-1/2} \right)}^{\left(x_{i+1/2},y_{j+1/2} \right)} \left(a_{11}\,\dfrac{\partial v }{\partial x}+ a_{12}\,\dfrac{\partial v }{\partial y}+ x\,b_1\,v \right)d y \\ \nonumber
    		&- \int_{\left(x_{i-1/2},y_{j-1/2} \right)}^{\left(x_{i-1/2},y_{j+1/2} \right)} \left(a_{11}\,\dfrac{\partial v }{\partial x}+a_{12}\,\dfrac{\partial v }{\partial y}+ x b_1\,v \right)d y \\ \nonumber
    		& + \int_{\left(x_{i-1/2},y_{j+1/2} \right)}^{\left(x_{i+1/2},y_{j+1/2} \right)} \left(a_{21}\,\dfrac{\partial v }{\partial x}+a_{22}\,\dfrac{\partial v }{\partial y} + y\, b_2\,v \right)d x\\ \nonumber
    		&- \int_{\left(x_{i-1/2},y_{j-1/2} \right)}^{\left(x_{i+1/2},y_{j-1/2} \right)} \left(a_{21}\,\dfrac{\partial v }{\partial x}+ a_{22}\,\dfrac{\partial v }{\partial y}+y\, b_2\,v \right)d x.
    	\end{align}
    	We shall look at (\ref{vol}) term by term. For the first term we want to approximate the integral by a constant as
    	\begin{align}
    		&\int_{\left(x_{i+1/2},y_{j-1/2} \right)}^{\left(x_{i+1/2},y_{j+1/2} \right)} \left(a_{11}\,\dfrac{\partial v }{\partial x}+ a_{12}\,\dfrac{\partial v }{\partial y}+ x \,b_1\,v \right)d y \\ \nonumber
    		& \approx \left(a_{11}\,\dfrac{\partial v }{\partial x}+ a_{12}\,\dfrac{\partial v}{\partial y}+x\, b_1\,v \right)|_{\left(x_{i+1/2},y_{j} \right)} \cdot h_{y_j}.
    	\end{align}
    	To achieve this, it is clear that we now need to derive approximations of the  $ k(v(x,y,t,\alpha_{i,j})) \cdot \bf n $ defined above at the mid-point $ \left(x_{i+1/2},y_{j} \right) $, of the interval $ I_{x_{i}} $ for $ i=0, 1,\cdots N_1-1$. This discussion is divided into two cases for $ i \geq 1 $ and $ i\,\in I_0 = (0,x_{1}) $. This is really an extension of the one dimensional fitted finite volume presented in the previous section.
	
	\textbf{\underline{Case I}:} For $ i\geq 2 $.
	
Remember that $ a_{11}(x,y,t,\alpha) = a(x, y, t, \alpha)\,x^2 $, we approximate the term $ \left(a_{11} \dfrac{\partial v}{\partial x}+ x\,b_1\,v \right) $ by solving the following two points boundary value problem
    	\begin{align}\label{vie2}	
		 \left({a} (x_{i+1/2}, y_j,t,\alpha_{i,j})\,x_{i+1/2} \dfrac{\partial v}{\partial x} + b_1 (x_{i+1/2},y_{j},t,\alpha_{i,j})\,v \right)'= 0 \\ \nonumber
    		\hspace*{2.5cm}  v(x_i,y_j,t)= v_{i,j}(t),\,\,\,\, v(x_{i+1},y_j,t)= v_{i+1,j}(t).
    	\end{align}
    	Integrating (\ref{vie2}) yields the first-order linear equations
    	\begin{align}
    		 {a} (x_{i+1/2},y_j,t,\alpha_{i,j})\,x_{i+1/2} \dfrac{\partial v}{\partial x}+ b_1 (x_{i+1/2},y_j,t,\alpha_{i,j})\,v  = C_1
    	\end{align}
    	where $ C_1 $ denotes an additive constant. Following the one dimensional fitted finite volume presented in the previous section, we have 
	\begin{align}
    		C_1 = \dfrac{{b_1}_{i+1/2,j}(t,\alpha_{i,j})\,\left(x_{i+1}^{\beta_{i,j}(t)}\,v_{i+1,j}-x_{i}^{\beta_{i,j}(t)}v_{i,j} \right)}{x_{i+1}^{\beta_{i,j}(t)}-x_{i}^{\beta_{i,j}(t)}}.
    	\end{align}
    	Therefore,
    	\begin{align}
    		a_{11}\,\dfrac{\partial v }{\partial x} + a_{12}\,\dfrac{\partial v }{\partial y}+ x\,b_1\,v \approx x_{i+1/2}\left( \dfrac{{b_1}_{i+1/2,j}(t,\alpha_{i,j})\,\left(x_{i+1}^{\beta_{i,j}(t)}v_{i+1,j}-x_{i}^{\beta_{i,j}(t)}\,v_{i,j} \right)}{x_{i+1}^{\beta_{i,j}(t)}-x_{i}^{\beta_{i,j}(t)}} + d_1\,y\,\,\dfrac{\partial v }{\partial y}\right),
    	\end{align}
where $ \beta_{i,j}(t)=\dfrac{{b_1}_{i+1/2,j}(t,\alpha_{i,j})}{{a}_{i+1/2,j}(t,\alpha_{i,j})} $ and $ a_{12} = a_{21} = d_{1}(x,y,t,\alpha)\,x\,y $.
	 Finally, we use the forward difference,
    	\begin{align*}
    		\dfrac{\partial v }{\partial y}\approx\dfrac{v_{i,j+1}-v_{i,j}}{h_{y_j}}
    	\end{align*}
    	Finally,
    	\begin{align}\label{vol1}
    		&\left[a_{11}\,\dfrac{\partial v }{\partial x}+ a_{12}\,\dfrac{\partial v}{\partial y}+ x\,b_1\,v \right]_{\left(x_{i+1/2},y_j \right)}
    		\cdot h_{y_j} \nonumber\\
    		&\approx x_{i+1/2}\left( \dfrac{{b_1}_{i+1/2,j}(t,\alpha_{i,j})\,\left(x_{i+1}^{\beta_{i,j}(t)}v_{i+1,j}-x_{i}^{\beta_{i,j}(t)}\,v_{i,j} \right)}{x_{i+1}^{\beta_{i,j}(t)}- x_{i}^{\beta_{i,j}(t)}}
    		+ {d_1}_{i,j}(t,\alpha_{i,j})\,y_j\,\dfrac{v_{i,j+1}-v_{i,j}}{h_{y_j}} \right) \cdot h_{y_j}.
    	\end{align}
Simillary, the second term in (\ref{vol}) can be approximated by
    	\begin{align}\label{vol2}
    		&\left[a_{11}\,\dfrac{\partial v }{\partial x}+ a_{12}\,\dfrac{\partial v}{\partial y}+x\,b_1\,v \right]_{\left(x_{i-1/2},{y_{j}} \right)}
    		\cdot h_{y_j}  \nonumber\\
    		&\approx x_{i-1/2}\left( \dfrac{{b_1}_{i-1/2,j}(t,\alpha_{i,j})\,\left(x_{i}^{\beta_{i-1,j}(t)}v_{i,j}-x_{i-1}^{\beta_{i-1,j}(t)}v_{i-1,j} \right)}{x_{i}^{\beta_{i-1,j}(t)}- x_{i-1}^{\beta_{i-1,j}(t)}}
    		+ {d_1}_{i,j}(t,\alpha_{i,j})\,y_j\,\dfrac{v_{i,j+1}-v_{i,j}}{h_{{y_j}}} \right) \cdot h_{y_j}.
    	\end{align}
    	
    	\textbf{\underline{Case II}:} For $ j\geq 2 $.\\
    	
    	For the third term we want to approximate the integral by a constant,  that is
    	\begin{align}
    		&\int_{\left(x_{i-1/2},y_{j+1/2} \right)}^{\left(x_{i+1/2},y_{j+1/2}\right)} \left(a_{21}\,\dfrac{\partial v }{\partial x}+ a_{22}\,\dfrac{\partial v }{\partial y}+ y\,b_2\,v \right)d x \\ \nonumber
    		& \approx \left(a_{21}\,\dfrac{\partial v }{\partial x}+ a_{22}\,\dfrac{\partial v}{\partial y}+ y\,b_2\,v \right)|_{\left({y_{j+1/2}},x_{i} \right)} \cdot h_{x_i}.
    	\end{align}
    	Following the first case of this section, we have 
 \begin{align}\label{vol3}
    		&\left[a_{21}\,\dfrac{\partial v }{\partial x}+ a_{22}\,\dfrac{\partial v}{\partial x_2}+ y\,b_2\,v \right]_{\left(x_{i},{y_{j+1/2}} \right)}
    		\cdot h_{x_i}  \nonumber\\
    		&\approx y_{j+1/2}\left( \dfrac{{b_2}_{i,j+1/2}(t,\alpha_{i,j})\,\left({y_{j+1}}^{\bar{\beta}_{i,j}(t)}v_{i,j+1}-{y_{j}}^{\bar{\beta}_{i,j} (t)}\,v_{i,j} \right)}{{y_{j+1}}^{\bar{\beta}_{i,j}(t)}- {y_{j}}^{\bar{\beta}_{i,j}(t)}}
    		+ {d_1}_{i,j}(t,\alpha_{i,j})\,x_i\,\dfrac{v_{i+1,j}-v_{i,j}}{h_{x_i}} \right) \cdot h_{x_i}.
    	\end{align}
     Similary, the fourth term in (\ref{vol}) can be approximated by
    	\begin{align}\label{vol4}
    		&\left[a_{21}\,\dfrac{\partial v }{\partial x}+ a_{22}\,\dfrac{\partial v}{\partial y}+y\,b_2\,v \right]_{\left(x_{i},y_{j-1/2} \right)}
    		\cdot h_{x_i} \approx  \nonumber\\
    		& y_{j-1/2}\left( \dfrac{{b_2}_{i,j-1/2}(t,\alpha_{i,j})\,\left(y_{j}^{\bar{\beta}_{i,j-1}(t)}\,v_{i,j}-y_{j-1}^{\bar{\beta}_{i,j-1}(t)}\,v_{i,j-1} \right)}{y_{j}^{\bar{\beta}_{i,j-1}(t)}- y_{j-1}^{\bar{\beta}_{i,j-1}(t)}}
    		+ {d_1}_{i,j}(t,\alpha_{i,j})\,x_i\,\dfrac{v_{i+1,j}-v_{i,j}}{h_{x_i}} \right) \cdot h_{x_i},
    	\end{align}
    	for $ j= 2,\cdots, N_2-1 $, where $ \bar{\beta}_{i,j}(t) = \dfrac{{b_2}_{i,j+1/2}(t,\alpha_{i,j})}{\bar{a}_{i,j+1/2}(t,\alpha_{i,j})} $ with $ a_{22}(x,y, t,\alpha) = \bar{a}(x,y,t,\alpha)\,y^2 $.  

    	\textbf{\underline{Case III:}} Approximation of the flux at $ I_0 $.
    	Note that the analysis in case I does not apply to the approximation of the flux on $[0,x_1] $ because (\ref{vie2}) is degenerated. Therefore, we reconsider the following form
    	\begin{align}\label{vie4}
    		(a(x_{1/2},y_j,t,\alpha_{1,j})\,x_{1/2} \dfrac{\partial v}{\partial x}+{b_1}(x_{1/2},y,t,\alpha_{1,j})\,v )' \equiv C_2\,\,\,\textbf{in}\,\,[0,x_{1}] \\ \nonumber
    		v(x_0,y_j)= v_{0,j},\,\,\,\, v(x_{1},{y_j})= v_{1,j},
    	\end{align}
    	where $ C_2 $ is an unknown constant to be determined. Integrating (\ref{vie4}), we can deduce that 
    	\begin{align}
    		&\left[a_{11}\,\dfrac{\partial v }{\partial x} + a_{12}\,\dfrac{\partial v}{\partial y}+ x\,b_1\,v \right]_{\left(x_{1/2},{y_{j}} \right)}
    		\cdot h_{y_j}  \nonumber\\
    		&\approx x_{1/2}\left( \dfrac{1}{2}\left[ (a_{x_{1/2},j}(t,\alpha_{1,j}) + {b_1}_{x_{1/2},j}(t,\alpha))\,v_{1,j}-(a_{x_{1/2},j}(t,\alpha_{1,j}) - {b_1}_{x_{1/2},j}(t,\alpha_{1,j}))\,v_{0,j}\right]\right. \\ \nonumber
    		&\left.+ {d_1}_{1,j}(t,\alpha_{1,j})\,y_j\,\dfrac{v_{1,j+1}-v_{1,j}}{h_{y_j}} \right) \cdot h_{y_j}.
    	\end{align}
    	\textbf{\underline{Case IV:}} Approximation of the flux at $ J_0 $.
    	
    	Using the same procedure for the approximation of the flux at $ I_0 $, we deduce that 
    	\begin{align}
    		&\left[a_{21}\,\dfrac{\partial v }{\partial x}+ a_{22}\,\dfrac{\partial v}{\partial y}+ y\,b_2\,v \right]_{\left(x_{i},{y_{1/2}} \right)}
    		\cdot h_{x_{i}}\approx  \nonumber\\
    			& y_{1/2}\left( \dfrac{1}{2}\left[(\bar{a}_{i,y_{1/2}}(t,\alpha_{i,1}) +{b_2}_{i,y_{1/2}}(t,\alpha))\,v_{i,1}-(\bar{a}_{i,y_{1/2}}(t,\alpha_{i,1})
    			-{b_2}_{i,y_{1/2}}(t,\alpha_{i,1}))\,v_{i,0}\right] \right. \\ \nonumber	&\left. +{d_1}_{i,1}(t,\alpha_{i,1})\,x_i\,\dfrac{v_{i+1,1}-v_{i,1}}{h_{x_i}} \right) \cdot h_{x_i}.
    	    	\end{align}
By replacing the flux by his value for $ i = 1,\cdots,N_1-1 $ and $ j = 1,\cdots,N_2-1 $, equation (\ref{bjr}) becomes
    	{\small{
    	\begin{align}
    		&\dfrac{d v_{i,j} }{d t} \\ \nonumber
    		& + \sup_{\alpha_{i,j} \in \mathcal{A}} \dfrac{1}{l_{i,j}}\,\left[x_{{i+1/2}}\left( \dfrac{{b_1}_{i+1/2,j}(t,\alpha)\,\left(x_{{i+1}}^{\beta_{i,j}(t)}\,v_{i+1,j}-x_{{i}}^{\beta_{i,j}(t)}\,v_{i,j} \right)}{x_{{i+1}}^{\beta_{i,j}(t)}- x_{{i}}^{\beta_{i,j}(t)}}
    		+ {d_1}_{i,j}(t,\alpha_{i,j})\,y_j\,\dfrac{v_{i,j+1}-v_{i,j}}{h_{y_j}} \right) \cdot h_{y_j} \right. \nonumber\\
    		& \left.- x_{{i-1/2}}\left( \dfrac{{b_1}_{i-1/2,j}(t,\alpha_{i,j})\,\left(x_{i}^{\beta_{i-1,j}}\,v_{i,j}-x_{{i-1}}^{\beta_{i-1,j}(t)}\,v_{i-1,j} \right)}{x_{{i}}^{\beta_{i-1,j}(t)}- x_{{i-1}}^{\beta_{i-1,j}(t)}}
    		+ {d_1}_{i,j}(t,\alpha_{i,j})\,y_j\,\dfrac{v_{i,j+1}-v_{i,j}}{h_{{y_j}}} \right) \cdot h_{{y_j}} \right.\nonumber\\
    		& \left. + {y_{j+1/2}}\left( \dfrac{{b_2}_{i,j+1/2}(t,\alpha_{i,j})\,\left({y_{j+1}}^{\bar{\beta}_{i,j}(t)}\,v_{i,j+1}-{y_{j}}^{\bar{\beta}_{i,j}(t)}v_{i,j} \right)}{{y_{j+1}}^{\bar{\beta}_{i,j}(t)} - {y_{j}}^{\bar{\beta}_{i,j}(t)}} + {d_1}_{i,j}(t,\alpha_{i,j})\,x_i\,\dfrac{v_{i+1,j}-v_{i,j}}{h_{x_{i}}} \right) \cdot h_{x_{i}} \right. \nonumber \\
    		& \left. -  {y_{j-1/2}}\left( \dfrac{{b_2}_{i,j-1/2}(t,\alpha_{i,j})\,\left({y_{j}}^{\bar{\beta}_{i,j-1}(t)}v_{i,j}-{y_{j-1}}^{\bar{\beta}_{i,j-1}(t)}\,v_{i,j-1} \right)}{{y_{j}}^{\bar{\beta}_{i,j-1}(t)}- {y_{j-1}}^{\bar{\beta}_{i,j-1}(t)}} \right.\right. \nonumber \\
    		&\left.\left. + {d_1}_{i,j}(t,\alpha_{i,j})\,x_i\,\dfrac{v_{i+1,j}-v_{i,j}}{h_{x_{i}}} \right) \cdot h_{x_{i}} + c_{i,j}(t,\alpha_{i,j})\, v_{i,j}\,l_{i,j}\right] = 0 \nonumber
    	\end{align}
    }}
    By setting $ \tau=T-t $ and including the boundary conditions, we have  the following  system 
    	\begin{align}\label{p5}
    		\begin{cases}
    			 & \underset{\alpha \in \mathcal{A}^N}{\sup}\, \left[e_{i-1,j}^{i,j}(\tau,\alpha)\, v_{i-1,j} +e_{i,j}^{i,j}(\tau,\alpha)\, v_{i,j} + e_{i+1,j}^{i,j}(\tau,\alpha)\, v_{i+1,j}+ e_{i,j-1}^{i,j} (\tau,\alpha)v_{i,j-1}+e_{i,j+1}^{i,j} (\tau,\alpha)v_{i,j+1} \right] \\
    			& -\dfrac{d v_{i,j} }{d \tau} = 0,\,\,\mbox{with}~~ \,\,\,\,  v(0) \,\,\,\text{given},
    		\end{cases}
    	\end{align}
	where  for $ i= 1,\cdots, N_1-1 $, $ j = 1,\cdots, N_2-1 $ and $N=(N_1-1)\times(N_2-1)$, we have
    	\begin{align}
    		e_{0,j}^{1,j}& = - \dfrac{1}{4\,l_{1,j} } h_{{y_j}}\,x_{1}(a_{x_{1/2},j}(\tau,\alpha_{1,j})- {b_1}_{x_{1/2},j}(\tau,\alpha_{1,j}))\,v_{0,j}\\ \nonumber
    		e_{1,j}^{1,j} & = \dfrac{1}{4\,l_{1,j} } h_{{y_j}}\,x_{1}(a_{x_{1/2},j}(\tau,\alpha_{1,j}) + {b_1}_{x_{1/2},j}(\tau,\alpha_{1,j})) - \dfrac{1}{2}\,c_{1,j}(\tau,\alpha_{1,j}) + {d_1}_{1,j}(\tau,\alpha_{1,j})\,x_i\dfrac{h_{{y_j}}}{l_{1,j}} \nonumber\\
    		& + x_{{1+1/2}} h_{{y_j}} \dfrac{{b_1}_{1+1/2,j}(\tau,\alpha_{1,j})\,x_{{1}}^{\beta_{1,j}(\tau)} }{l_{1,j}\left(x_{{2}}^{\beta_{1,j}(\tau)}- x_{{1}}^{\beta_{1,j}(\tau)}\right)} \\
    		e_{2,j}^{1,j} & = -{d_1}_{1,j}(\tau,\alpha_{1,j})\,x_i\,\dfrac{h_{{y_j}}}{l_{1,j}} - x_{{1+1/2}} h_{{y_j}} \dfrac{{b_1}_{1+1/2,j}(\tau,\alpha_{1,j})\,x_{{2}}^{\beta_{1,j}(\tau)} }{l_{1,j}\left(x_{{2}}^{\beta_{1,j}(\tau)}- x_{{1}}^{\beta_{1,j}(\tau)}\right)}\\
    		e_{i,0}^{i,1}& = -\dfrac{1}{4\,l_{i,1} } h_{x_{i}}\,y_{1}( \bar{a}_{i,y_{1/2}}(\tau,\alpha_{i,1}) - {b_2}_{i,y_{1/2}}(\tau,\alpha_{i,1})) \,v_{i,0}\\ \nonumber
    		e_{i,1}^{i,1}& = \dfrac{1}{4\,l_{i,1} } h_{x_{i}}\,y_{1}( \bar{a}_{i,y_{1/2}}(\tau,\alpha_{i,1}) +{b_2}_{i,y_{1/2}}(\tau,\alpha_{i,1}) - \dfrac{1}{2}\,c_{i,1}(\tau,\alpha_{i,1}) + {d_1}_{i,1} (\tau,\alpha_{i,1})\,y_j\,\dfrac{h_{x_{i}}}{l_{i,1}}  \nonumber \\
    		& + y_{{1+1/2}} h_{x_{i}} \dfrac{{b_2}_{i,1+1/2}(\tau,\alpha_{i,1})\,y_{{1}}^{\bar{\beta}_{i,1}(\tau)} }{l_{i,1}\left(y_{{2}}^{\bar{\beta}_{i,1}(\tau)}- y_{{1}}^{\bar{\beta}_{i,1}(\tau)}\right)} \\
    		e_{i,2}^{i,1}& = -{d_1}_{i,1}(\tau,\alpha_{i,1}) \,y_j\,\dfrac{h_{x_{i}}}{l_{i,1}} - {y_{1+1/2}} h_{x_{i}} \dfrac{{b_2}_{i,1+1/2}(\tau,\alpha_{i,1})\,{y_{2}}^{\bar{\beta}_{i,1}(\tau)} }{l_{i,1}\left({y_{2}}^{\bar{\beta}_{i,1}(\tau)}- {y_{1}}^{\bar{\beta}_{i,1}(\tau)}\right)}\\
    		e_{i+1,j}^{i,j} &= - {d_1}_{i,j}(\tau,\alpha_{i,j}))\,\dfrac{h_{{y_j}}}{l_{i,j}} - x_{{i+1/2}} h_{{y_j}} \dfrac{{b_1}_{i+1/2}(\tau,\alpha_{i,j})\,x_{{i+1}}^{\beta_{i,j}(\tau)} }{l_{i,j}\left(x_{{i+1}}^{\beta_{i,j}(\tau)}- x_{{i}}^{\beta_{i,j}(\tau)}\right)}\\
    		e_{i-1,j}^{i,j}& = - x_{{i-1/2}} h_{{y_j}} \dfrac{{b_1}_{{i-1/2,j}}(\tau,\alpha_{i,j}))\,x_{{i-1}}^{\beta_{i-1,j}(\tau)} }{l_{i,j}\left(x_{{i}}^{\beta_{i-1,j}(\tau)}- x_{{i-1}}^{\beta_{i-1,j}(\tau)}\right)}\\
    		e_{i,j}^{i,j}& = {d_1}_{i,j}(\tau, \alpha_{i,j}))\,x_i\,\dfrac{h_{{y_j}}}{l_{i,j}} + x_{{i+1/2}} h_{{y_j}} \dfrac{{b_1}_{{i+1/2,j}}(\tau,\alpha_{i,j}))\,x_{{i}}^{\beta_{i,j}(\tau)} }{l_{i,j}\left(x_{{i+1}}^{\beta_{i,j}(\tau)}- x_{{i}}^{\beta_{i,j}(\tau)}\right)} \\ \nonumber
    	  &	+ x_{{i-1/2}} h_{{y_j}} \dfrac{{b_1}_{{i-1/2,j}}(\tau,\alpha_{i,j}))\,x_{{i}}^{\beta_{i-1,j}(\tau)} }{l_{i,j}\left(x_{{i}}^{\beta_{i-1,j}(\tau)}- x_{{i-1}}^{\beta_{i-1,j}(\tau)}\right)}-c_{i,j}(\tau,\alpha_{i,j})) \\ \nonumber
    	  & {d_1}_{i,j}(\tau,\alpha_{i,j})) \,y_j\,\dfrac{h_{x_{i}}}{l_{i,j}} + {y_{j+1/2}} h_{x_{i}} \dfrac{{b_2}_{i,j+1/2}(\tau,\alpha_{i,j}))\,{y_{j}}^{\bar{\beta}_{i,j}(\tau)} }{l_{i,j}\left({y_{j+1}}^{\bar{\beta}_{i,j}(\tau)}- {y_{j}}^{\bar{\beta}_{i,j}(\tau)}\right)} \\ \nonumber
    	   & +{y_{j-1/2}} h_{x_{i}} \dfrac{{b_2}_{i,j-1/2}(\tau,\alpha_{i,j}))\,{y_{j}}^{{\beta}_{i,j-1}(\tau)} }{l_{i,j}\left({y_{j}}^{\bar{\beta}_{i,j-1}(\tau)}- {y_{j-1}}^{\bar{\beta}_{i,j-1}(\tau)}\right)}  
    		\end{align}
    \begin{align}
		e_{i,j+1}^{i,j} &= - {d_1}_{i,j}(\tau,\alpha_{i,j})) \,y_j\,\dfrac{h_{x_{i}}}{l_{i,j}} - y_{{j+1/2}} h_{x_{i}} \dfrac{{b_2}_{{i,j+1/2}}(\tau,\alpha_{i,j}))\,{y_{j+1}}^{\bar{\beta}_{i,j}(\tau)} }{l_{i,j}\left({y_{j+1}}^{\bar{\beta}_{i,j}(\tau)}- {y_{j}}^{\bar{\beta}_{i,j}(\tau)}\right)}\\
    		e_{i,j-1}^{i,j}& = - {y_{j-1/2}} h_{x_{i}} \dfrac{{b_2}_{i,j-1/2}(\tau,\alpha_{i,j}))\,{y_{j-1}}^{\bar{\beta}_{i,j-1}(\tau)} }{l_{i,j}\left({y_{j}}^{\bar{\beta}_{i,j-1}(\tau)}- {y_{j-1}}^{\bar{\beta}_{i,j-1}(\tau)}\right)}.
    	\end{align}
	As for one dimension case, \eqref{p5} can be rewritten as the Ordinary Differential Equation (ODE) coupled with optimization problem
    	\begin{align}\label{p4}
 		\begin{cases}
   		 \dfrac{d \textbf{v}(\tau)}{d \tau} + \underset{\alpha \in \mathcal{A}^N}{\inf}\,\left[E (\tau, \alpha)\,\textbf{v}(\tau) + F(\tau,,\alpha) \right] = 0, \\
   			~~~\mbox{with}~~ \,\,\,\,  \textbf{v}(0)\,\,\,\text{given},
   		\end{cases}
    	\end{align} or 
    	\begin{align}\label{pan1}
    		\begin{cases}
    			\dfrac{d \textbf{v}(\tau)}{d \tau}  =\underset{\alpha \in \mathcal{A}^N}{\sup}\,\left[A (\tau,\alpha)\,\textbf{v}(\tau) + G(\tau, \alpha) \right] \\
    			~~~\mbox{with}~~ \,\,\,\,  \textbf{v}(0) \,\,\,\text{given},
    		\end{cases}
    	\end{align}
    	where  $ A(\tau,\alpha) = - E (\tau,\alpha) $,
	 $ \textbf{v} = \left({v}_{1,1},\cdots,{v}_{1,N_2-1},\cdots ,{v}_{N_1-1,1},\cdots,{v}_{N_1-1,N_2-1}\right) $ and  $G (\tau, \alpha) =-F(\tau,\alpha)$ includes boundary condition.
	 \begin{thm}
    \label{tm2} 
    		Assumme that $ c_{i,j}(\tau,\alpha) < 0 $, \,\, ${d_1}_{i,j}(\tau,\alpha) > 0,\; i = 1,\cdots, N_1-1 $,\,\, $ j = 1,\cdots, N_2-1 $, and let  $h =\underset{\underset{1\leq j\leq N_2}{1\leq i\leq N_1}} {\max} l_{i,j}$.
		 If $h$ is reatively small then the matrix $E (\tau,\alpha) =\left(e_{i,j}^{i,j}\right)_{\underset{j = 1,\cdots ,N_2-1,}{i = 1,\cdots ,N_1-1}} $ in (\ref{p4}) is an \textbf{M}-matrix for any $ \alpha\, \,\in\,\mathcal{A}^N$. 
    	\end{thm}
    	\hspace*{0.5cm}
	\begin{proof} 
	The Proof follows the same lines of that of \thmref{tm2}.
		\end{proof}
	\section{Temporal discretization and optimization }  
	\label{sec2}
  	This section is devoted to the numerical time discretization method  for the spatially discretized optimization problem using the fitted finite volume method. 
	We will present it in one and two dimensional.
Let us re-consider the differential equation coupled with optimization problem given in \eqref{p4xx} or \eqref{pan1} by
\begin{align}
\label{nmod}
\textbf{v}_{\tau}(\tau) = \sup_{\alpha \in \mathcal{A}^N} \left[ A (\tau,\alpha) \textbf{v}(\tau)+ G(\tau,\alpha) \right]\\ \nonumber
~~\,\,\,\textbf{v}(0)\,\,\,\text{given},
\end{align}
For temporal discretization, we use a constant time step $\Delta t > 0$, of course variable time steps can be used.  
The temporal grid points given by~ $\Delta t = \tau_{n+1}-\tau_n $~ for ~ $n =  1, 2,\ldots m-1 $. We denote $\textbf{v}(\tau_n) \approx \textbf{v}^n$, $A^n (\alpha)=A(\tau_n,\alpha)$ and $G^n (\alpha)=G(\tau_n,\alpha).$ 

For $\theta \,\in \left[\dfrac{1}{2}, 1\right]$, following \cite{HPFH}, the $\theta$-Method approximation of \eqref{nmod} in time is given by
\begin{align}\label{scheme}
& \textbf{v}^{n+1} - \textbf{v}^{n} = \Delta t \,\sup_{\alpha \in \mathcal{A}^N} \left( \theta\,  [A^{n+1} (\alpha)\,\textbf{v}^{n+1} + G^{n+1}(\alpha)] \right. \\ \nonumber
&\left. + (1-\theta)\, [A^{n} (\alpha)\,\textbf{v}^{n} + G^{n}(\alpha)]\right).
\end{align}
As we can notice, to find the unknown $\textbf{v}^{n+1}$ we need also to solve an optimization.  Let 
\begin{align}
\label{opt}
\alpha^{n+1} \in \left(\underset{\alpha \in \mathcal{A}^N }{arg \sup} \left\lbrace \theta \,\Delta t \left[ A^{n+1}(\alpha)\,\textbf{v}^{n+1} +  G^{n+1}(\alpha )\right] + (1-\theta)\,\Delta t \left[A^{n} (\alpha )\,\textbf{v}^{n} + G^{n}(\alpha)\right] \right\rbrace\right).
\end{align}
Then, the unknown $\textbf{v}^{n+1}$ is solution  of the following equation
\begin{align}\label{vu1}
& [ I  - \theta\, \Delta t \,A^{n+1} (\alpha^{n+1})]\,\textbf{v}^{n+1} = [I + (1-\theta)\,\Delta t\, A^{n} (\alpha^{n+1})]\,\textbf{v}^{n} \\
&+[\theta\, \Delta t \, G^{n+1}(\alpha^{n+1})+(1-\theta) \Delta t\, G^{n}(\alpha^{n+1})] \nonumber,
\end{align} 
Note that, for $\theta= \dfrac{1}{2}$, we have the \textit{Crank Nickolson scheme} and for $\theta=1$ we have the \textit{Implicit scheme}.
Unfortunately \eqref{scheme}-\eqref{opt} are nonlinear and coupled and we need to iterate at every time step.
The following  iterative scheme close to the one in \cite{HPFH} is used.
\begin{enumerate}
	\item Let  $ \left( \textbf{v}^{n+1}\right)^0=\textbf{v}^{n}$, 
		\item Let $ \hat{\textbf{v}}^{k}= \left( \textbf{v}^{n+1}\right)^k$,
	\item  For $k=0,1,2 \cdots $ until convergence ($\Vert \hat{\textbf{v}}^{k+1}-\hat{\textbf{v}}^{k}\Vert \leq \epsilon$, given tolerance) solve
	\begin{align} \label{mi2}
	&\alpha^{k}_i \in \left(\underset{\alpha \in \mathcal{A}^N }{arg \sup} \left\lbrace \theta \,\Delta t \left[ A^{n+1} (\alpha)\,\hat{\textbf{v}}^k+  G^{n+1}(\alpha) \right]_i  + (1-\theta)\,\Delta t \, \left[A^{n} (\alpha )\,\textbf{v}^{n} + G^{n}(\alpha)\right]_i \right\rbrace\right)&\\
	 & \alpha^{k}=(\alpha^{k})_i&\\
	& [ I  - \theta\, \Delta t\,A^{n+1} (\alpha^{k})]\,\hat{\textbf{v}}^{k+1} = [I + (1-\theta)\,\Delta t \, A ^{n}(\alpha^{k})] \textbf{v}^{n} \\
	&+[\theta\, \Delta t \, G^{n+1}(\alpha^{k})+(1-\theta) \Delta t \, G^{n}(\alpha^{k})] \nonumber,
	\end{align}
	\item  Let  $k_l$ being the last iteration in step 3,  set $\textbf{v}^{n+1}:=\hat{\textbf{v}}^{k_l}$,\,\, $\alpha^{n+1}:=\alpha^{k_l}$.

\end{enumerate}	
\section{Numerical experiments}
\label{sec3}
The goal of this section is   carried out on test problems in both 1 and 2 space dimensions to validate the numerical scheme presented in the previous section.
  All computations were performed in Matlab 2013  using the estimate parameters coming from \cite{JH} and \cite{HPFH}.  
  We will present two problems with exact solution and one problem without exact solution modelling  cash management in finance.
\begin{prob}
	Consider the following Merton's stochastic control problem  such that  $ \alpha = \alpha(t,x)$ is a feedback control belongs in $ (0,1) $ 
		
		\begin{align} \label{pb1}
		\begin{cases}
		v(x,t)=\underset{\alpha \in\, (0,1)}{\max} \mathbb{E}\left\lbrace \dfrac{1}{p}\,x^p(T)\right\rbrace.\\
		dx_t =  b^{\alpha_t}(t,x_t)\, dt + \sigma^{\alpha_t} (t,x_t)\,d\omega_t\\ 
		\end{cases}	
		\end{align}
		where $ b^{\alpha_t}(t,x_t) = 	x_t\,\left[ r + \alpha_t\, (\mu-r)\right] $, $\sigma^{\alpha_t} (t,x_t)= x_t\, \alpha_t \,\sigma  $,
		$  0 < p < 1 $ is coefficient of risk aversion, $r$, $\mu$, $\sigma $ are constants, $ x_t\,\in\,\mathbb{R} $ and $ \omega_t $  Brownian motion. We assume  $ \mu > r $. For this problem, the corresponding HJB equation is given by
		\begin{align}\label{pb2}
		\begin{cases}
		v_t(t, x) + \underset{\alpha \in (0,1)}{\sup} \left[L^{\alpha} v(t, x)\right] = 0 \quad\text{on} \  [0,T)\times \mathbb{R}\\
		v(T,x) = \dfrac{x^p}{p}, \,\, \,x \,\in \mathbb{R}
		\end{cases}
		\end{align} 
		where 
		\begin{align}
		L^{\alpha} v(t, x) =  (b^{\alpha}(t, x)) \dfrac{\partial v(t,x)}{\partial x}  +  ( a^{\alpha} (t,x)) \,\dfrac{\partial^2 v(t,x)}{\partial x^2},
		\end{align}   
		and $ a^\alpha (t,x) = \dfrac{1}{2}(\sigma^{\alpha}(t,x))^2 $.
		
		The divergence form of the HJB (\ref{pb2}) is given  by
	\begin{align}
	\dfrac{\partial v (t, x)}{\partial t} + \sup_{\alpha \in \,(0,1)}\left[ \dfrac{\partial }{\partial x} \left(a (t,x,\alpha)\,x^2 \dfrac{\partial v(t, x)}{\partial x}+ b (t,x,\alpha)\,x\,v(t, x) \right)+c(t, x,\alpha)\,v(t, x)\right] = 0,
	\end{align}
where 
	\begin{align*}
	a(t, x, \alpha) &= \dfrac{1}{2}\sigma^2\alpha^2\\
	b(t, x,\alpha)&= r+ (\mu-r)\,\alpha-\sigma^2\alpha^2\\
	c(t, x, \alpha)&=-(r+\alpha\, (\mu-r)-\sigma^2\,\alpha^2).
	\end{align*}
		
\end{prob}
 The Domain where we compare the solution is $ \Omega =\left[ 0, x_{{max}}\right]$,  where Dirichlet Boundary conditions is used at the boundaries.
  Of course the value of the boundary  conditions are taken to be equal to the exact solution.  The exact  solution given in \cite{HP} is given at every $(x_i, \tau^n)$  by
\begin{eqnarray}
v\left(\tau^n, x_i\right) &=&  e^{p \times (n \times \Delta t- T) \times \rho} \times \dfrac{(x_i)^p}{p},\\
\rho &=& \left( r+\dfrac{(\mu-r)^2}{\sigma^2\,(1-p)} + \dfrac{1}{2}(p-1)\,\sigma^2\,\left(\dfrac{(\mu-r)}{\sigma^2\,(1-p)}\right)^2  \right),\,\,\,\,0 < p < 1
\end{eqnarray}
We use the  following $L^2(\Omega \times [0,T])$ norm of the absolute error 
\begin{eqnarray}
\ \left\| v^m-v\right\|_{L^2[\Omega \times [0,T]]} = \left(\sum_{n = 0}^{m-1} \sum_{i = 1}^{N_1-1}  (\tau_{n+1}-\tau_n)\times l_{i}\times (v_{i}^n - v\left( \tau^n, x_{i}, \right))^2 \right)^{1/2},
\end{eqnarray}
where $v^m$ is the numerical approximation of $v$ computed from  our numerical scheme.

\begin{figure}[!h]
        			\begin{minipage}[b]{0.5\linewidth}
        				\centering\includegraphics[scale=0.5]{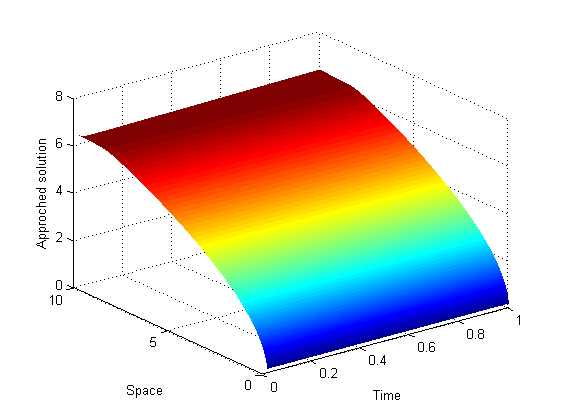}
        				\caption{Implicit Fitted Finite Volume.}
        				\label{eq 261b}
        			\end{minipage}\hfill
        			\begin{minipage}[b]{0.5\linewidth}
        				\centering\includegraphics[scale=0.5]{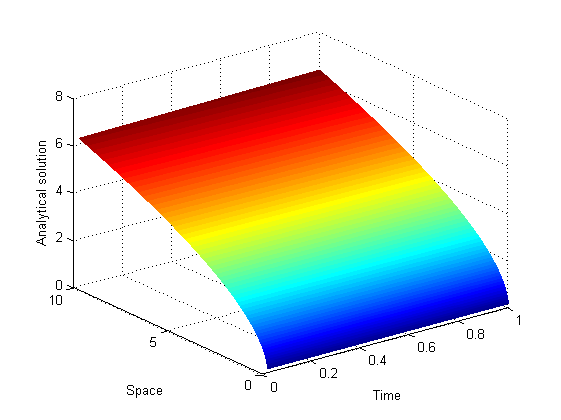}
        				\caption{Ansatz Analytical solution }
        				\label{eq 281b}
        			\end{minipage}
        		\end{figure}

  The 3 D  graphs of the Implicit Fitted Finite Volume ( $\theta=1$)  with its corresponding exact solution is given at  Figure \ref{eq 261b} and  Figure \ref{eq 281b}.
  For our computation, we have  $[0,10]$  for computational domain  with  $N= 1500$, $r = 0.0449$, $\mu = 0.0657$, $\sigma = 0.2537$, $p= 0.5255$ and   $T=1 $.  

\begin{figure}[!h]
	\begin{minipage}[b]{0.5\linewidth}
		\centering\includegraphics[scale=0.5]{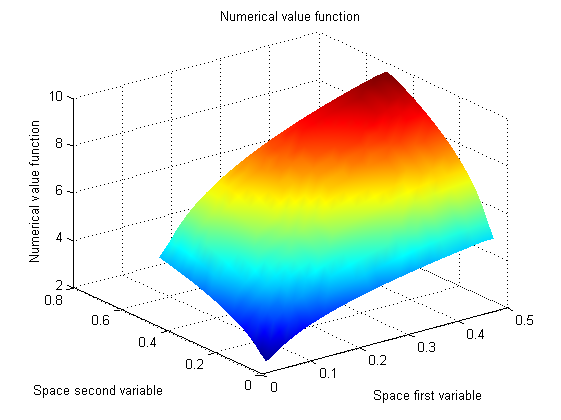}
		\caption{Implicit Fitted Finite Volume ($\theta=1$) at finite time $T=1$.}
		\label{eq261b}
	\end{minipage}\hfill
	\begin{minipage}[b]{0.5\linewidth}
		\centering\includegraphics[scale=0.5]{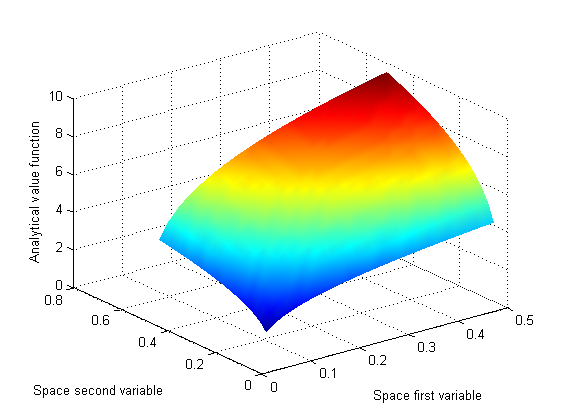}
		\caption{Ansatz Analytical solution at finite time $T=1$ }
		\label{eq281}
	\end{minipage}
\end{figure}

   We compare  the fitted finite volume method and the finite difference method in Table \ref{2i}
   \newpage
   \vspace{1cm}
     \begin{table}[!h]
   	\begin{center}

   \begin{tabular}{|c|c|c|c|c|}
   	\hline
   	Time subdivision &  $200 $ &  $150 $ & $ 100 $ &  $50 $ \\ 
   	\hline Error of Implicit Fitted Finite Volume method & 3.34 E-01  & 6.81 E-01  & 1.01 E-00 & 1.33 E-00\\ 
   	 \hline Error of Implicit Finite difference method & 3.37 E-01  & 6.89 E-01  & 1.02 E-00 & 1.34 E-00 \\  
   	 \hline
   	 \end{tabular}
   	 \vspace*{1cm}
   \caption{Comparison of the  implicit fitted finite volume method and implicit finite difference method.}
   \label{2i}
\end{center}
\end{table}
From  Table \ref{2i}, we can observe the accuracy of the  implicit fitted finite volume comparing to the implicit finite difference method.
\begin{prob}
	Consider the  following two  dimensional Merton's stochastic control model  such that  $ \alpha_1 = \alpha_1(t,x)$ and $ \alpha_2 = \alpha_2(t,y)$ are a feedback control in $ (0,1) $ 
	\begin{align} \label{pb11}
	v(t,x,y)=\underset{(\alpha_1, \alpha_2) \in\, (0,1)\times (0,1)}{\max} \mathbb{E}\left\lbrace \dfrac{1}{p}\,x^p(T) \times \dfrac{1}{p}\,y^p(T)\right\rbrace,
	\end{align}
	\begin{align}
	\mbox{subject}\,\,\text{ to}\,\, \begin{cases}
	dx_t =  {b_1}^{{\alpha_1}_t}(t,x_t)\, dt + \sigma^{{\alpha_2}_t} (t,x_t)\,d{\omega_1}_t\\
	dy_t =  {b_2}^{{\alpha_2}_t}(t,y_t)\, dt + \sigma^{{\alpha_2}_t} (t,y_t)\,d{\omega_2}_t\\ 
	\end{cases}	
	\end{align}
	where 
	\begin{align*}
&	{b_1}^{{\alpha_1}_t}(t,x_t) = 	x_t\,\left[ r_1 + {\alpha_1}_t\, (\mu_1-r_1)\right],\\
& {b_2}^{{\alpha_2}_t}(t,y_t) = 	y_t\,\left[ r_2 + {\alpha_2}_t\, (\mu_2-r_2)\right],\\
&\sigma^{{\alpha_1}_t} (t,x_t)= x_t\, {\alpha_1}_t \,\sigma,\,\quad \quad  \sigma^{{\alpha_2}_t} (t,y_t)= y_t\, {\alpha_2}_t \,\sigma,
	\end{align*} 
	$  0 < p < 1 $ is coefficient of risk aversion, $r_1$, $\mu_1$, $r_2$, $\mu_2$ $\sigma $ are constants, $ x_t\, y_t\,\in\,\mathbb{R} $ and $ \rho $ the correlation of the two Brownian motion. We assume that $ \mu_1 > r_1 $ and $ \mu_2 > r_2 $. For this problem, the corresponding HJB equation is given by
	\begin{align}\label{pb6}
	\begin{cases}
	v_t(t, x, y) + \underset{(\alpha_1, \alpha_2) \in (0, 1)\times(0, 1 )}{\sup} \left[L^{\alpha_1,\alpha_2} v(t, x, y)\right] = 0 \quad\text{on} \  [0,T)\times \mathbb{R}\times \mathbb{R}\\
	v(T, x , y) = \dfrac{x^p}{p}\times \dfrac{y^p}{p}, \,\, \,x,\, \,y\,\in \mathbb{R}
	\end{cases}
	\end{align} 
	where 
\begin{align*}
& L^{\alpha} v(t, x, y) =  ({b_1}^{\alpha_1}(t, x)) \dfrac{\partial v(t,x,y)}{\partial x} + ({b_2}^{\alpha_2}(t, y)) \dfrac{\partial v(t,x, y)}{\partial y}  + \dfrac{1}{2}(\sigma^{\alpha_1}(t,x))^2\,\dfrac{\partial^2 v(t,x ,y )}{\partial x^2}\\
& + \dfrac{1}{2}(\sigma^{\alpha_2}(t, y))^2\,\dfrac{\partial^2 v(t,x ,y )}{\partial y^2} + (\sigma^{\alpha_1}(t,x))\,(\sigma^{\alpha_2}(t, y))\,\dfrac{\partial^2 v(t,x ,y )}{\partial x^2},
\end{align*}   
and the two dimensional divergence form is given by 
\begin{align} 
\dfrac{\partial v(t, x, y ) }{\partial t} + \sup_{\alpha \in \,(0,1)\times(0,1)}\left[ \nabla\cdot \left( k(t, x, y,\alpha)\right) + c(t,x,y,\alpha)\,v(t,x,y) \right] = 0,
\end{align}
where\,\,\,
$k(t, x, y,\alpha) = A(t, x, y,\alpha)\cdot\nabla v(t, x, y)+ b(t,x,y,\alpha)\cdot v(t,x,y)$ 
\[
A = \left( \begin{array}{cc}
a_{11} & a_{12} \\
a_{21} & a_{22}
\end{array} \right),
\]
\begin{eqnarray*}
a_{11}(t,x,y,\alpha) = \dfrac{1}{2}\sigma^2\,\alpha_1^2 \,x^2 ,~ a_{22}(t,x,y,\alpha) = \dfrac{1}{2} \,\sigma^2\,\alpha_2^2 \,y^2, ~ a_{12}(t,x,y,\alpha) = a_{21}(x,y,t,\alpha) = \dfrac{1}{2}  \sigma^2\,\alpha_1 \,\alpha_2 \,x\,y.
\end{eqnarray*}
By identification,
\begin{eqnarray*}
a(t,x,y,\alpha) &=& \dfrac{1}{2}\sigma^2\,\alpha_1^2\\
\bar{a}(t,x,y,\alpha)&=& \dfrac{1}{2} \,\sigma^2\,\alpha_2^2\\
b_1(t,x,y,\alpha) &=& r_1 + \alpha_1\, (\mu_1-r_1) - \dfrac{1}{2}\sigma^2\,\alpha_1\,\alpha_2- \sigma^2\,\alpha_1^2, \\
b_2(t,x,y,\alpha) & =&  r_2 + \alpha_2\, (\mu_2-r_2) -\dfrac{1}{2}\sigma^2\,\alpha_1\,\alpha_2-\sigma^2\,\alpha_2^2,\\
c (t,x,y,\alpha)& =& -\left[ r_1 + (\mu_1-r_1)\,\alpha_1\right]- \left[ r_2 + (\mu_2-r_2)\,\alpha_2 \right] +\sigma^2\left(\alpha_1^2  +\alpha_2^2 + \alpha_1 \,\alpha_2\right),\\
d_1(t,x,y,\alpha)&=&\dfrac{1}{2}  \sigma^2\,\alpha_1 \alpha_2.\nonumber
\end{eqnarray*}
The two dimensional  Antsaz  exact solution \cite{HP}  at  $\left(\tau^n, x_i, y_j\right) $ is  given by
{\small{
\begin{eqnarray*}
&&v\left(\tau^n, x_i, y_j\right) =  e^{p \times (n \times \Delta t - T) \times \rho} \times \dfrac{(x_i)^p}{p} \times \dfrac{(y_j)^p}{p},\\
&&\rho = \underset{\alpha_1,\alpha_2 \,\in (0,1)\times(0,1)}{\sup}\left[ r_1 + r_2 +  (\mu_1-r_1)\,\alpha_1 + (\mu_2-r_2)\,\alpha_2+ \dfrac{1}{2} \sigma^2\,\alpha_1^2\,(p-1) 
 + \dfrac{1}{2} \sigma^2\,\alpha_2^2\,(p-1)
+ \sigma^2\,\alpha_1\,\alpha_2\,p \right],\\
&& 0 < p < 1/2.
\end{eqnarray*}
}}
\end{prob}
We use the  following $L^2(\Omega \times [0,T]),\,\, \Omega = [0,x_{{\text{max}}}] \times [0,y_{{\text{max}}}]$ norm of the absolute error 
\begin{eqnarray}
\ \left\| v^m-v\right\|_{L^2[\Omega \times [0,T]]} = \left(\sum_{n = 0}^{m-1} \sum_{i = 1}^{N_1-1} \sum_{j = 1}^{N_2-1}  (\tau_{n+1}-\tau_n)\times hx_{i}\times hy_{j}\times (v_{i,j}^n - v\left( \tau^n, x_{i}, x_{j},\right))^2 \right)^{1/2},
\end{eqnarray}
where $v^m$ is the numerical approximation of $v$ computed from  our numerical scheme.
%
The 3 D  graphs of the Implicit Fitted Finite Volume ( $\theta=1$ at  the final time $T=1$)  with its corresponding exact solution is given at  Figure \ref{eq261b} and  Figure \ref{eq 281b}, with 
$N_1= 50$, $N_2 = 45$, $r_1 = 0.0449/2$, $\mu_1 = 0.0657/2$, $r_2 = 0.044/2$, $\mu_2 = 0.065/2$, $\sigma = 0.2537/2$ and $p= 0.5255/2$.
%
 We compare  the fitted finite volume method and the finite difference method in Table \ref{1i}.
 Again, we can observe the accuracy of the fitted scheme comparing to the finite difference scheme.
\begin{table}[!h]
	\begin{center}
		
		\begin{tabular}{|c|c|c|c|c|}
			\hline
			Time subdivision &  $200 $ &  $150 $ & $ 100 $ &  $50 $ \\ 
			\hline Error of Fitted Finite Volume method & 4.08 E-02  & 7.84 E-02  & 1.14 E-01 & 1.47 E-01\\ 
			\hline Error of Finite difference method & 4.23 E-02  & 7.93 E-02  & 1.16 E-01 & 1.48 E-01\\   
			\hline
		\end{tabular}
		\vspace*{1cm}
		\caption{Errors  table  for  fitted finite volume method and finite difference method in dimension 2.}
		\label{1i}
	\end{center}
\end{table}
%

 	
 	\begin{prob}
 	We consider a optimal Cash Management under a stochastic volatility Model problem coming from \cite{kylan}. We assume that the firm invests its cash in a bank account and a stock in a portfolio of
 	value $w_t$ at time $t$, and the proportion of wealth invested in the stock at time $t$ is $u_t$. The interest rate earned in the bank account is $ r_1 $, the return from the stock at time
 	$ t $ has two components, the cash dividend rate $ r_2  $, the capital gain rate $ R_t $. The dynamic of the capital gain rate $ R_t $ is assumed to be governed by the stochastic process
 	\begin{align}
 	dR_t &= [\beta_1\, R_t + f]\, dt +  \sigma_t  \,d{W_{1}}_t,
 	\end{align}
 	 and the volatility $ \sigma_t $ with modeled by
 	 \begin{align}
 	 	d\sigma_t &= \alpha\,\sigma_t\, dt + \beta\,\sigma_t \, d{W_{2}}_t. 
 	 \end{align} 
 	 Suppose that the firm has a demand rate $ d_t $ for cash at time $ t $, and that the demand
 	rate $ d(t) $ is normally distributed with mean $ 0 $ and variance $ 0.2 $. We assume that  $ u_t\,\in\,[0,1] $ and 
 	the
 	wealth dynamics for this cash management problem is given by
\begin{align}
dw_t &= w_t \,u_t\,r_2\,dt +w_t\,(1-u_t)\,r_1\,dt + w_t\,R_t\,dt -d(t)\,w_t\,dt.
\end{align}
The objective of the firm is to maximize the expectation of the total holdings at the
terminal time $ T $. The portfolio optimization problem is given by
 		\begin{align} \label{pbr2}
 		J(w, R, \sigma, T) =\underset{u \in [0,1]}{\max} \mathbb{E}\left\lbrace w_T\right\rbrace.
 		\end{align}
 		\begin{align}
 			\mbox{subject}\,\,\text{ to}\,\, \begin{cases}
 			dw_t &= w_t \,u_t\,r_2\,dt +w_t\,(1-u_t)\,r_1\,dt + w_t\,R_t\,dt -d(t)\,w_t\,dt,\,\,\\ \nonumber
 			dR_t &= [\beta_1\, R_t + f]\, dt +  \sigma_t  \,d{W_{1}}_t\\ \nonumber
 				d\sigma_t &= \alpha\,\sigma_t\, dt + \beta\,\sigma_t \, d{W_{2}}_t \nonumber
 			\end{cases}	
 		\end{align}
 		We assume that the two Brownian motions are correlated,  that is  $ d{W_{1}}_t\,d{W_{2}}_t = \rho\,dt $.
 		 For this problem of  optimal Cash Management  the analytical solution is not available, so our numerical scheme will  to provide approximated solution. 
		The  corresponding HJB equation for this  optimal cash management problem  (\ref{pbr2}) is given by 
 		\begin{align}\label{pbr1}
 		& J_t +  \underset{u \in [0,1]}{\max}  \left\lbrace (f + \beta_1\,R)\,J_R +   w\,(u\,r_2 + (1-u)\,r_1 + u\,R - d(t) )\,J_w +  \right.\\
 		& \left. 1/2\,\left(\sigma^2\,J_{RR}  +\beta^2\,\sigma^2 \,J_{\sigma \sigma} + 2\,\rho\,\beta\,\sigma^2\,J_{\sigma\,R} \right)  \nonumber
 		+ \alpha\,\sigma\,J_\sigma   \right\rbrace =0 ,
 		\end{align}
 		with terminal condition  $J( \cdot, T) = w_T. $ 
		To simplify the problem, we assume that $$J(w, R, \sigma, t) = w H( R, \sigma, t).$$   Therefore (\ref{pbr1}) is equivalent to
 		\begin{align}\label{84}
 			& H_t +  \underset{u \in [0,1]}{\max}  \left\lbrace (f + \beta_1\,R)\,H_R +   (u\,r_2 + (1-u)\,r_1 + u\,R - d (t))\,H +  \right.\\
 			& \left. 1/2\,\left(\sigma^2\,J_{RR}  +\beta^2\,\sigma^2 \,H_{\sigma \sigma} + 2\,\rho\,\beta\,\sigma^2\,H_{\sigma\,R} \right)  \nonumber
 			+ (\alpha\,\sigma)\,H_\sigma   \right\rbrace =0 
 		\end{align}
 		with terminal condition  $H( R, \sigma, T) = 1 $. The HJB
 		equation (\ref{84}) is a problem with two state variables  $R $ and  $\sigma$. The divergence form of the problem (\ref{84}) is then  given by
 			\begin{align} 
 			\dfrac{\partial H(R,\sigma,t) }{\partial t} + \sup_{u \,\in \,[0,1]}\left[ \nabla\cdot \left( k(R,\sigma,t,u)\right) + c(R,\sigma,t,u)\,H(R,\sigma,t) \right] = 0,
 			\end{align}
 			where\,\,\,
 			$k(R,\sigma,t,u) = A(R,\sigma,t,u)\cdot\nabla H(R,\sigma,t)+ b(R,\sigma,t,u)\cdot H(R,\sigma,t)$ 
 			\[
 			A=\left( \begin{array}{cc}
 			a_{11} & a_{12} \\
 			a_{21} & a_{22}
 			\end{array} \right),
 			\]
 			\begin{align*}
 			a_{11} = \dfrac{1}{2}\sigma^2 ,~ a_{22} = \dfrac{1}{2} \,\beta^2\, \,\sigma^2, ~ a_{12} = a_{21} = \dfrac{1}{2}\,  \sigma\,\rho\, \beta.
 			\end{align*}
 			By identification,
 			\begin{align*}
 			a(R,\sigma,t) &= \dfrac{\sigma^2}{2\,R^2},\,\,\,
 			\bar{a}(R,\sigma,t)= \dfrac{\beta^2}{2}\\
 			b_1(R,\sigma,t) &= \dfrac{f}{R} + \beta_1 - \dfrac{ \rho\,\beta}{2\,R},\,\,\,\,
 			b_2(R,\sigma,t)  = \alpha -\beta^2  \\
 			c (R,\sigma,t,u)& = u\,r_2 +(1-u)\,r_1 +u\,R -d(t) -\beta_1-\alpha +\beta^2.\nonumber
 			\end{align*}			
 	 	\end{prob} 
 	 	Because we have a stochastic volatility model, to solve  the PDE equation, we  have  considered the  following boundary conditions of Heston model
 	 	 \begin{align}
 	 	 H(0, \sigma, t) = 0,\\
 	 	 H(R, \sigma_{\text{max}}, t) = R,\\
 	 	 \dfrac{\partial H }{\partial R} (R_{\text{max}}, \sigma, t) = 1.
 	 	 \end{align}
 	 	Because the PDE has two second derivatives in  the two spatial directions, four boundary
 	 	conditions are needed. This comes from the fact that the two second order derivatives give rise
 	 	to two unknown integration constants. To meet this requirement, at the boundary  $\sigma = 0 $ it is
 	 	considered inserting $\sigma = 0 $ into the PDE to complete the set of four boundary conditions:
 	 	\begin{align}
 	 	& H_t(R, 0, t) +  \underset{u \in [0,1]}{\max}  \left\lbrace (f + \beta_1\,R)\,H_R(R,0,t) +   (u\,r_2 + (1-u)\,r_1 + u\,R - d )\,H(R,0,t)   \right\rbrace =0 
 	 	 	\end{align}
 	 	The HJB equation is
 	 solved using some parameters  values  in \cite{kylan}  given in  the following  tabular 
 	 
%
 \hspace*{1cm} 
\begin{center}
 	  	 \begin{tabular}{|c|c|c|c|c|c|c|}
 	 	\hline $ f $ & $ \beta_1 $  & $ \beta $ & $ \alpha $  & $ r_1 $ & $ r_2 $ & $ \rho $  \\	 	 
 	 	\hline
 	 	 $ 0.12 $ & $ 0.96 $ & $ 0.3 $ & $ -0.85 $ &  $0.024 $ & $ 0.01 $ & $ 0.5 $  \\
 	 	\hline 
\end{tabular} 
\vspace*{2cm}	\\

\end{center}
 	 \begin{figure}[!h]
 	 		\centering\includegraphics[scale=1]{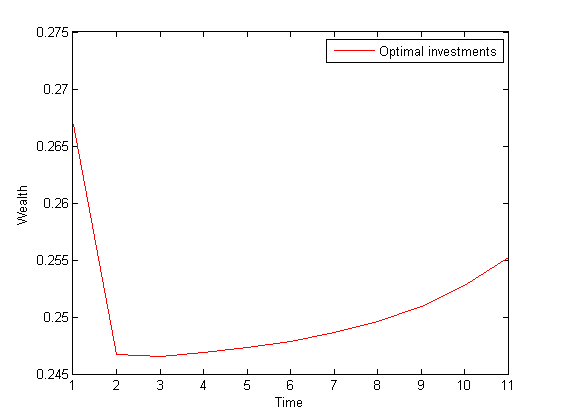}
 	 		\caption{  A fitted finite volume sample   solution of the wealth rate $H$ at  the point  $1/2,1/2$}.
 	 		\label{a1}
 	 \end{figure}
	 
 	 Figure \ref{a1} shows  a sample of  fitted finite volume  solution  of the wealth rate $H$  at the point $(1/2,1/2)$ from  $t = 1$  to $ t = 10 $ with $ N_1 = 10 $,  $N_2 = 10 $, $ R_{\text{max}} = 1/2 $, $ \sigma_{\text{max}} = 1/2 $. We can estimate the mean and moment of $H$   using Monte Carlo Method by generating many samples of $H$.
 	 \section{Conclusion}
	 \label{sec4}
 	 We presented a fitted  finite volume method  to solve  the HJB equation from stochastic optimal control problems  coupled with implicit temporal discretization.  The optimization  problem  is solved  at every time step using iterative method.
	 It was shown that the corresponding  system matrix is an M-matrix, so the maximum principle is preserved for the discrete system.
	 Numerical experiments in 1 and 2 dimensions are performed to prove the accuracy of the fitted finite volume method comparing to the standard finite difference methods.
 \section*{Acknowledgements}
 The first author  was supported by the project African Center of Excellence in Mathematics and Applied Sciences (ACE-MSA) in Benin and the European Mathematical Society (EMS).

\end{document}